\documentclass{amsart}%
\usepackage{amsfonts}
\usepackage{amsmath}
\usepackage{amssymb}
\usepackage{graphicx}%
\newtheorem{theorem}{Theorem}
\theoremstyle{plain}

\newtheorem{corollary}[theorem]{Corollary}

\newtheorem{lemma}[theorem]{Lemma}

\newtheorem{remark}[theorem]{Remark}
\begin{document}
\title[The existence of suitable vector fields]{A remark on the existence of suitable vector fields related to the dynamics of
scalar semi-linear parabolic equations}
\author{Fengbo Hang}
\address{Department of Mathematics, Michigan State University, East Lansing, MI 48824}
\email{fhang@math.msu.edu}
\author{Huiqiang Jiang}
\address{School of Mathematics, University of Minnesota, 127 Vincent Hall, 206 Church
St. S.E., Minneapolis, MN 55455}
\email{hqjiang@math.umn.edu}
\thanks{The research of the first author is supported in part by NSF Grant DMS-0209504.}
\date{March 25, 2005}

\begin{abstract}
In 1992, P. Pol\'{a}\v{c}ik\cite{P2} showed that one could
linearly imbed any vector fields into a scalar semi-linear
parabolic equation on $\Omega$ with Neumann boundary condition
provided that there exists a smooth vector field $\Phi=\left(
\phi_{1},\cdots,\phi_{n}\right)  $ on $\overline{\Omega}$ such
that
\[
\left\{
\begin{array}
[c]{l}%
\operatorname*{rank}\left(  \Phi\left(  x\right)  ,\partial_{1}\Phi\left(
x\right)  ,\cdots,\partial_{n}\Phi\left(  x\right)  \right)  =n\text{ for all
}x\in\overline{\Omega},\\
\frac{\partial\Phi}{\partial\nu}=0\text{ on }\partial\Omega\text{.}%
\end{array}
\right.
\]
In this short note, we give a classification of all the domains on which one
may find such type of vector fields.

\end{abstract}
\maketitle

Let $\Omega$ be a bounded smooth domain in $\mathbb{R}^{n}$, $\nu$ be the
outer normal direction of $\partial\Omega$ and $f\in C^{\infty}\left(
\overline{\Omega}\times\mathbb{R\times R}^{n},\mathbb{R}\right)  $. The
infinite dimensional dynamical system defined by%
\begin{equation}
\left\{
\begin{tabular}
[c]{l}%
$u_{t}=\Delta u+f\left(  x,u,\nabla u\right)  ,\quad x\in\Omega,t>0,$\\
$\frac{\partial u}{\partial\nu}=0$ on $\partial\Omega,$%
\end{tabular}
\right. \label{system}%
\end{equation}
on suitable Sobolev spaces has attracted much interest (see
\cite{A,H} and more recent references at the end of this note).
When $n=1$, equation $(\ref{system})$ has rather simple dynamics,
each bounded solution will converge to an equilibrium. The
situation is quite different when $n\geq 2$, the solutions of
$(\ref{system})$ can exhibit very complicated behavior. One
relatively easy way to demonstrate the complexity of its dynamical
behavior is the realization of ODEs in $(\ref{system})$. We refer
the reader to the recent survey paper\cite{P4} by P.
Pol\'{a}\v{c}ik for a quick overview of the progress made up to
2002, more results can be found in
\cite{DP,P1,P2,P3,P5,PR,Pr,PrR1,PrR2,PrR3,R} and the references
therein.

In particular, the following nice result was proved in \cite{P2}:
if there exists a smooth vector field $\Phi$ on
$\overline{\Omega}$, $\Phi=\left(  \phi
_{1},\cdots,\phi_{n}\right)  $ such that
\[
\left\{
\begin{array}
[c]{l}%
\operatorname*{rank}\left(  \Phi\left(  x\right)  ,\partial_{1}\Phi\left(
x\right)  ,\cdots,\partial_{n}\Phi\left(  x\right)  \right)  =n\text{ for all
}x\in\overline{\Omega},\\
\frac{\partial\Phi}{\partial\nu}=0\text{ on }\partial\Omega\text{,}%
\end{array}
\right.
\]
then for any smooth vector field $X$ on $\mathbb{R}^{n}$, there exists a
smooth function $f$, such that the linear space $\operatorname*{span}\left\{  \phi_{1}%
,\cdots,\phi_{n}\right\}  $ is invariant under (\ref{system}) and for any
integral curve of $X$, $c=c\left(  t\right)  $, $u=\sum_{i=1}^{n}c_{i}\left(
t\right)  \phi_{i}\left(  x\right)  $ is a solution to (\ref{system}).
Moreover, it was shown that such kind of vector field always exists on a
starshaped domain. The main result of this short note is a classification of
all the domains on which one may find this type of vector fields. More
precisely, we have

\begin{theorem}
\label{thm1}Let $\Omega\subset\mathbb{R}^{n}$ $\left(  n\geq2\right)  $ be an
open bounded smooth domain, then the necessary and sufficient condition for
the existence of a smooth map $F:\overline{\Omega}\rightarrow\mathbb{R}^{n}$
with
\begin{equation}
\left\{
\begin{array}
[c]{l}%
\operatorname*{rank}\left(  F\left(  x\right)  ,\partial_{1}F\left(  x\right)
,\cdots,\partial_{n}F\left(  x\right)  \right)  =n\text{ for any }%
x\in\overline{\Omega},\\
\frac{\partial F}{\partial\nu}=0\text{ on }\partial\Omega\text{,}%
\end{array}
\right. \label{eqmain}%
\end{equation}
is that $\overline{\Omega}$ is diffeomorphic to $\overline{B}_{1}$
or $\overline{B}_{2}\backslash B_{1}$. Here $B_{1}$ and $B_2$ are
two open balls centered at zero with radius $1$ and $2$
respectively.
\end{theorem}

\begin{remark}
In fact, if $\overline{\Omega}$ is diffeomorphic to
$\overline{B}_{1}$, then any solution to $(\ref{eqmain})$, $F$,
must have exactly one zero in $\Omega$. If $\overline{\Omega}$ is
diffeomorphic to $\overline{B}_{2}\backslash B_{1}$, then any
solution to $(\ref{eqmain})$, $F$, does not vanish at all. These
conclusions will follow from the arguments below.
\end{remark}
\begin{remark}
Our theorem improves the linear imbedding result of P.
Pol\'{a}\v{c}ik \cite{P2}, and at the same time, it exhibits the
limitation of linear realization method. On the other hand,
realization of ODEs in a nonlinear fashion can be made in a much
more general setting. For example, it is shown in
\cite{P3,PrR2,PrR3} that any ODE, in any dimension, has an
arbitrarily small perturbation that is realizable nonlinearly in a
single semi-linear parabolic equation defined on any given open
subset of $\mathbb{R}^n$, $n\geq 2$.
\end{remark}
To prove our main result, we first show that Neumann boundary
condition in $(\ref{eqmain})$ can be relaxed.

\begin{lemma}
\label{lem1}Let $\Omega\subset\mathbb{R}^{n}$ $\left(
n\geq2\right)  $ be an open bounded smooth domain, if there exists
a smooth map $G:\overline{\Omega }\rightarrow\mathbb{R}^{n}$ such
that
\[
\operatorname*{rank}\left(  G\left(  x\right)  ,\partial_{1}G\left(  x\right)
,\cdots,\partial_{n}G\left(  x\right)  \right)  =n\text{ for any }%
x\in\overline{\Omega}%
\]
and%
\[
\dim\operatorname*{span}\left\{  G\left(  x\right)  ,\operatorname*{im}\left(
\left.  G\right\vert _{\partial\Omega}\right)  _{\ast,x}\right\}  =n\text{ for
any }x\in\partial\Omega,
\]
here $\left(  \left.  G\right\vert _{\partial\Omega}\right)  _{\ast,x}$
denotes the tangent map of $\left.  G\right\vert _{\partial\Omega}$ at $x$,
then we may find a smooth map $F:\overline{\Omega}\rightarrow\mathbb{R}^{n} $
satisfying (\ref{eqmain}).

\begin{proof}
Let $\varepsilon>0$ be small enough such that the map%
$$
\phi :\partial\Omega\times\left[  0,3\varepsilon\right]
\rightarrow \left\{  y\in\overline{\Omega}:\text{dist}\left(
y,\partial\Omega\right) \leq3\varepsilon\right\} $$ defined by
$$\phi\left( x,t\right)= x-t\nu\left(  x\right)
$$
is a diffeomorphism. Let $P=G\circ\phi$. For any $t\in\left[  0,3\varepsilon
\right]  $, let $P_{t}\left(  x\right)  =P\left(  x,t\right)  $ for
$x\in\partial\Omega$, then we may assume $\varepsilon$ is small enough such
that for any $t\in\left[  0,3\varepsilon\right]  $,%
\[
\dim\operatorname*{span}\left\{  P_{t}\left(  x\right)  ,\operatorname*{im}%
\left(  P_{t}\right)  _{\ast,x}\right\}  =n\text{ for all }x\in\partial\Omega.
\]
Let $\eta:\mathbb{R\rightarrow R}$ be a smooth function such that%
\[
\eta\left(  t\right)  =\left\{
\begin{tabular}
[c]{l}%
$\varepsilon,$ when $t\leq\varepsilon/2;$\\
$t$, when $t\geq3\varepsilon/2;$%
\end{tabular}
\right.
\]
and $\eta^{\prime}\left(  t\right)  \geq0$ for all $t$. Define%
\[
Q_{t}\left(  x\right)  =Q\left(  x,t\right)  =P\left(  x,\eta\left(  t\right)
\right)  \text{ for }x\in\partial\Omega\text{, }0\leq t\leq3\varepsilon
\text{.}%
\]
Then it is clear that for any $t\in\left[  0,3\varepsilon\right]  $,%
\[
\dim\operatorname*{span}\left\{  Q_{t}\left(  x\right)  ,\operatorname*{im}%
\left(  Q_{t}\right)  _{\ast,x}\right\}  =n\text{ for any }x\in\partial\Omega.
\]
Let
\[
F\left(  y\right)  =\left\{
\begin{array}
[c]{lcl}%
G\left(  y\right)   & \text{if} & y\in\overline{\Omega},\, \operatorname*{dist}%
\left(  y,\partial\Omega\right)  \geq2\varepsilon,\\
Q\left(  \phi^{-1}\left(  y\right)  \right)   & \text{if} & y\in
\overline{\Omega},\, \operatorname*{dist}\left(
y,\partial\Omega\right) \leq3\varepsilon,
\end{array}
\right.
\]
then it is easy to check that $F$ satisfies all the requirements.
\end{proof}
\end{lemma}

Since the existence of a smooth vector field $G$ in Lemma 1 is a
property of the $\overline{\Omega}$ which is preserved under
diffeomorphisms, we conclude

\begin{corollary}
\label{cor1}Assume $\overline{\Omega}_{1}$ is diffeomorphic to $\overline
{\Omega}_{2}$, and for $\Omega_{1}$ we may find a solution to (\ref{eqmain}),
then we may find a solution to (\ref{eqmain}) for $\Omega_{2}$ too.
\end{corollary}

To derive the necessary condition for the existence of a vector field
satisfying (\ref{eqmain}), we will need

\begin{lemma}
\label{lem2}Let $\Omega\subset\mathbb{R}^{n}$ $\left(  n\geq2\right)  $ be an
open bounded smooth domain, if there exists a smooth map $H:\overline{\Omega
}\rightarrow S^{n-1}$ such that%
\[
\operatorname*{rank}\left(  \partial_{1}H\left(  x\right)  ,\cdots
,\partial_{n}H\left(  x\right)  \right)  =n-1\text{ for any }x\in
\overline{\Omega}%
\]
and%
\[
\dim\operatorname*{im}\left(  \left.  H\right\vert _{\partial\Omega}\right)
_{\ast,x}=n-1\text{ for any }x\in\partial\Omega,
\]
then $\overline{\Omega}$ is diffeomorphic to $\overline{B}_{2}\backslash
B_{1}$.

\begin{proof}
First we claim that each path connected component of $\partial\Omega$ is
diffeomorphic to $S^{n-1}$. This is clear when $n=2$. If $n\geq3$, since
$\left.  H\right\vert _{\partial\Omega}$ has full rank everywhere and
$\partial\Omega$ is compact, $$\left.  H\right\vert _{\partial\Omega}%
:\partial\Omega\rightarrow S^{n-1}$$ is a covering map (see
\cite{M}). Since $S^{n-1}$ is simply connected, we see each path
connected components of $\partial\Omega$ must be diffeomorphic to
$S^{n-1}$. Indeed, the restriction of $H$ to such a component
serves as a diffeomorphism.

To proceed, we observe that from the assumption on $H$, it follows from
implicit function theorem that for any $\xi\in S^{n-1}$, $H^{-1}\left(
\xi\right)  $ is a smooth one dimensional submanifold of $\overline{\Omega}$,
moreover $H:\overline{\Omega}\rightarrow S^{n-1}$ is a smooth fiber bundle
(see \cite{M}). Fix a point $x_{0}\in\partial\Omega$, let $\xi_{0}=H\left(
x_{0}\right)  $ and $\Gamma=H^{-1}\left(  \xi_{0}\right)  $, then we have an
exact sequence (see Theorem 6.7 of chapter VII in \cite{B})%
\[
\pi_{n-1}\left(  \Gamma,x_{0}\right)  \rightarrow\pi_{n-1}\left(
\overline{\Omega},x_{0}\right)  \rightarrow\pi_{n-1}\left(  S^{n-1},\xi
_{0}\right)  \rightarrow\pi_{n-2}\left(  \Gamma,x_{0}\right)  .
\]
If $n\geq3$, then both $\pi_{n-1}\left(  \Gamma,x_{0}\right)  $
and $\pi _{n-2}\left(  \Gamma,x_{0}\right)  $ vanishes, this shows
$$\pi_{n-1}\left( \overline{\Omega},x_{0}\right)
\widetilde{=}\mathbb{Z}$$ and hence $\overline{\Omega}$ is
diffeomorphic to $\overline{B}_{2}\backslash B_{1}$. If $n=2$,
then since $\pi_{1}\left(  \Gamma,x_{0}\right)  $ vanishes and
$\pi _{0}\left(  \Gamma,x_{0}\right)  $ is finite, we see
$\pi_{1}\left( \overline{\Omega},x_{0}\right)  $ is again
isomorphic to $\mathbb{Z}$, this shows $\overline{\Omega}$ must be
diffeomorphic to $\overline{B}_{2}\backslash B_{1}$.
\end{proof}
\end{lemma}

Now we are ready to prove the main theorem.

\begin{proof}
[Proof of theorem \ref{thm1}]First if $\Omega=B_{1}$ or
$B_{2}\backslash \overline{B}_{1}$, then $G\left(  x\right)  =x$
satisfies the assumption in the Lemma \ref{lem1}, half of the
theorem follows from Corollary \ref{cor1}. On the other hand,
assume for some $\Omega$, we may find a smooth map $F$ satisfying
(\ref{eqmain}). For $x\in\partial\Omega$, choose a base for
the tangent space of $\partial\Omega$ at $x$, namely $e_{1},\cdots,e_{n-1}$, then%

\begin{align*}
&  \operatorname*{rank}\left(  F\left(  x\right)  ,\partial_{1}F\left(
x\right)  ,\cdots,\partial_{n}F\left(  x\right)  \right) \\
 =& \operatorname*{rank}\left(  F\left(  x\right)  ,F_{\ast}e_{1}%
,\cdots,F_{\ast}e_{n-1},F_{\ast}\nu\right) \\
 =& \operatorname*{rank}\left(  F\left(  x\right)  ,F_{\ast}e_{1}%
,\cdots,F_{\ast}e_{n-1}\right)  =n,
\end{align*}
hence $F\left(  x\right)  \neq0$ for any $x\in\partial\Omega$.
Moreover, it follows from the fact
\[
\operatorname*{rank}\left(  F\left(  x\right)  ,\partial_{1}F\left(  x\right)
,\cdots,\partial_{n}F\left(  x\right)  \right)  =n\text{ for any }%
x\in\overline{\Omega}%
\]
that the zeroes of $F$ in $\Omega$ must be isolated. Hence $F$ has
at most finitely many zeroes, say $x_{1},\cdots,x_{m}$, here
$m\geq 0$. For $\varepsilon>0$ small enough, let
\[
U=\Omega\backslash%
{\bigcup\limits_{i=1}^{m}} \overline{B_{\varepsilon}\left(
x_{i}\right)  },
\]
if $m\geq 1$ and $U=\Omega$ if $m=0$, then we have
\[
\operatorname*{rank}\left(  F\left(  x\right)  ,\partial_{1}F\left(  x\right)
,\cdots,\partial_{n}F\left(  x\right)  \right)  =n\text{ for any }%
x\in\overline{U}%
\]
and%
\[
\dim\operatorname*{span}\left\{  F\left(  x\right)  ,\operatorname*{im}\left(
\left.  F\right\vert _{\partial U}\right)  _{\ast,x}\right\}  =n\text{ for any
}x\in\partial U.
\]
Let
\[
H\left(  x\right)  =\frac{F\left(  x\right)  }{\left\vert F\left(  x\right)
\right\vert }\text{ for }x\in\overline{U}\text{,}%
\]
then clearly%
\[
\operatorname*{rank}\left(  \partial_{1}H\left(  x\right)  ,\cdots
,\partial_{n}H\left(  x\right)  \right)  =n-1\text{ for any }x\in\overline{U}%
\]
and%
\[
\dim\operatorname*{im}\left(  \left.  H\right\vert _{\partial U}\right)
_{\ast,x}=n-1\text{ for any }x\in\partial U.
\]
It follows from Lemma \ref{lem2} that $\overline{U}$ must be
diffeomorphic to $\overline{B}_{2}\backslash B_{1}$, hence $m\leq
1$ and $\overline{\Omega} $ must be diffeomorphic to either
$\overline{B}_{1}$ or $\overline{B}_{2}\backslash B_{1}$.
\end{proof}

{\bf Acknowledgment}: We would like to thank Professor Wei-Ming Ni
for bringing the problem to our attention. We also thank Professor
Peter Pol\'{a}\v{c}ik for his kind comments on the manuscript.

\end{document}